\documentclass{gtart}
\usepackage{amsmath,amssymb}

\input gtmonout
\volumenumber{2}
\volumeyear{1999}
\volumename{Proceedings of the Kirbyfest}
\pagenumbers{103}{111}
\papernumber{5}
\received{11 April 1998}\revised{16 October 1999}
\published{17 November 1999}

\theoremstyle{plain}
\newtheorem{thm}{Theorem}[section]

\newtheorem{prop}[thm]{Proposition}

\newtheorem*{ther}{Theorem}

\theoremstyle{definition}

\def \R {\mathbf{R}}
\def \Z {\mathbf{Z}}
\def \Sig{\Sigma}

\def \vt {\vartheta}

\def \a {\alpha}
\def \b {\beta}
\def \g {\gamma}

\def \lam {\lambda}

\def \s {\sigma}
\def \t {\tau}

\def \bd {\partial}
\def \x {\times}
\def \ve {\varepsilon}

\def \ssw {\text{SW}}
\def \sw {\mathcal{SW}}

\def \DD {\Delta}

\def \f {\hat{f}}
\def \S {section }

\begin{document}

\title{Nondiffeomorphic Symplectic
$4$--Manifolds\\ with the same Seiberg--Witten Invariants}
\asciititle{Nondiffeomorphic Symplectic
4-Manifolds with the same Seiberg-Witten Invariants}

\shorttitle{Nondiffeomorphic Symplectic 4-Manifolds}

\authors{Ronald Fintushel\\Ronald J Stern}
\address{Department of Mathematics, Michigan State University\\
East Lansing, Michigan 48824, USA\\\smallskip\\Department of Mathematics,
University of California\\
Irvine,  California 92697, USA}
\asciiaddress{Department of Mathematics, Michigan State University,
East Lansing, Michigan 48824, USA\\Department of Mathematics,
University of California,
Irvine,  California 92697, USA}
\email{ronfint@math.msu.edu, rstern@math.uci.edu}

\begin{abstract}
The goal of this paper is to demonstrate that, at least for nonsimply
connected $4$--manifolds,
the Seiberg--Witten invariant alone does not determine diffeomorphism type
within the same homeomorphism type.
\end{abstract}

\asciiabstract{The goal of this paper is to demonstrate that, at least
for nonsimply connected 4-manifolds, the Seiberg-Witten invariant
alone does not determine diffeomorphism type within the same
homeomorphism type.}


\maketitle

\cl{\small\em Dedicated to Robion C Kirby on the
occasion of his
$60^{th}$ birthday}

\section{Introduction\label{Intro}}

The goal of this paper is to demonstrate that, at least for nonsimply
connected $4$--manifolds,
the Seiberg--Witten invariant alone does not determine diffeomorphism type
within the same homeomorphism type. The
first examples
which demonstrate this phenomenon were constructed by Shuguang Wang
\cite{ShW}. These are
examples of two homeomorphic $4$--manifolds with $\pi_1=\Z_2$ and trivial
Seiberg--Witten invariants. One of these
manifolds is irreducible and the other splits as a connected sum.
It is our goal here to exhibit examples among symplectic 4--manifolds, where
the Seiberg--Witten
invariants are known to be nontrivial. We shall construct symplectic
$4$--manifolds with
$\pi_1=\Z_p$ which have the same nontrivial Seiberg--Witten invariant but whose
universal covers have
different Seiberg--Witten invariants. Thus, at the very least, in order to
determine diffeomorphism type, one needs to consider the Seiberg--Witten
invariants of finite covers.

Recall that the Seiberg--Witten
invariant  of a smooth closed oriented $4$--manifold
$X$ with $b_2 ^+(X)>1$ is an integer-valued function which is defined on
the set of $spin ^{\, c}$ structures over $X$ (cf \cite{W}). In case
$H_1(X,\Z)$
has no 2--torsion  there
is a
natural identification of the $spin ^{\, c}$ structures of
$X$ with the characteristic elements of $H_2(X,\Z)$ (ie, those elements
$k$ whose
Poincar\'e duals $\hat{k}$ reduce mod~2 to $w_2(X)$).
In this case we view the
Seiberg--Witten invariant as
\[ \ssw_X\co \lbrace k\in H_2(X,\Z)|\hat{k}\equiv w_2(TX)\pmod2)\rbrace
\rightarrow \Z. \]
The sign of $\sw_X$
depends on an orientation of
$H^0(X,\R)\otimes\det H_+^2(X,\R)\otimes \det H^1(X,\R)$. If $\ssw_X(\b)\neq
0$, then $\b$ is called a {\it{basic class}} of $X$. It is a fundamental
fact that the set of
basic classes is
finite. Furthermore, if $\b$ is a basic class, then so is $-\b$ with
$\ssw_X(-\b)=(-1)^{(\text{e}+\text{sign})(X)/4}\,\ssw_X(\b)$ where
$\text{e}(X)$ is
the Euler number and $\text{sign}(X)$ is the signature of $X$.

Now let
$\{\pm \b_1,\dots,\pm \b_n\}$ be the set of nonzero basic classes for $X$.
Consider variables $t_{\b}=\exp(\b)$ for each $\b\in H^2(X;\Z)$ which
satisfy the
relations $t_{\a+\b}=t_{\a}t_{\b}$. We may then view the Seiberg--Witten
invariant
of $X$ as the Laurent polynomial
\[\sw_X = \ssw_X(0)+\sum_{j=1}^n \ssw_X(\b_j)\cdot
(t_{\b_j}+(-1)^{(\text{e}+\text{sign})(X)/4}\, t_{\b_j}^{-1}).\]
\medskip

\section{The Knot and Link Surgery  Construction}

We shall need the knot surgery construction of \cite{KL4M}: Suppose that we
are given a smooth
simply connected oriented $4$--manifold $X$ with $b^+> 1$ containing an
essential smoothly
embedded torus $T$ of self-intersection 0. Suppose further that
$\pi_1(X\setminus T)=1$ and
that $T$ is contained in a cusp neighborhood. Let $K\subset S^3$ be a
smooth knot and $M_K$ the
3--manifold obtained from 0--framed surgery on $K$. The meridional loop $m$
to $K$ defines a
1--dimensional homology class $[m]$ both in $S^3\setminus K$ and in $M_K$.
Denote by $T_m$
the torus $S^1\x m\subset S^1\x M_K$. Then $X_K$ is defined to be the fiber
sum
\[ X_K= X\#_{T=T_m} S^1\x M_K = (X\setminus N(T))\cup(S^1\x(S^3\setminus
N(K)), \]
where $N(T)\cong D^2\x T^2$ is a tubular neighborhood of $T$ in $X$ and
$N(K)$ is a
neighborhood of $K$ in $S^3$. If $\lam$ denotes the longitude of $K$
($\lam$ bounds a
surface in $S^3\setminus K$) then the gluing of this fiber sum identifies
$\{ {\text{pt}}
\}\x \lam$ with a normal circle to $T$ in $X$. The main theorem of
\cite{KL4M} is:
\begin{ther}{\rm\cite{KL4M}}\qua With the assumptions above, $X_K$ is
homeomorphic to $X$, and
\[ \sw_{X_K}=\sw_X\cdot \DD_K(t)\]
where $\DD_K$ is the symmetrized Alexander polynomial of $K$ and
$t=\exp(2[T])$.
\end{ther}

In case the knot $K$ is fibered, the $3$--manifold $M_K$ is a surface bundle
over the circle;
hence $S^1\x M_K$ is a surface bundle over $T^2$. It follows from \cite{Th}
that
$S^1\x M_K$ admits a symplectic structure and $T_m$ is a symplectic
submanifold. Hence, if
$T\subset X$ is a torus satisfying the conditions above, and if in addition
$X$ is a symplectic
$4$--manifold and $T$ is a symplectic submanifold,
then the fiber sum  $X_K= X\#_{T=T_m} S^1\x M_K$
carries a symplectic structure \cite{Gompf}. Since $K$ is a
fibered knot, its Alexander polynomial is the characteristic polynomial of
its monodromy
$\varphi$; in particular,
$M_K=S^1\x_{\varphi}\Sig$ for some surface $\Sig$ and
$\DD_K(t)=\det(\varphi_*-tI)$, where
$\varphi_*$ is the induced map on $H_1$.

There is a generalization of the above theorem in this case due to Ionel
and Parker \cite{IP} and to
Lorek \cite{L}.

\begin{ther}{\rm\cite{IP,L}}\qua Let $X$ be a symplectic $4$--manifold with
$b^+>1$, and let $T$ be a
symplectic self-intersection $0$ torus in $X$ which is contained in a cusp
neighborhood. Also,
let $\Sig$ be a symplectic 2--manifold with a symplectomorphism
$\varphi\co\Sig\to\Sig$ which has
a fixed point $\varphi(x_0)=x_0$. Let $m_0=S^1\x_{\varphi}\{ x_0\}$ and
$T_0=S^1\x m_0\subset
S^1\x(S^1\x_{\varphi}\Sig)$. Then
$X_{\varphi}=X\#_{T=T_0}S^1\x(S^1\x_{\varphi}\Sig)$ is a
symplectic manifold whose Seiberg--Witten invariant is
\[ \sw_{X_{\varphi}}=\sw_X\cdot \DD(t) \]
where $t=\exp(2[T])$ and $\DD(t)$ is the obvious symmetrization of
$\det(\varphi_*-tI)$.
\end{ther}
\noindent Note that in case $K$ is a fibered knot and
$M_K=S^1\x_{\varphi}\Sig$, Moser's theorem
\cite{M} guarantees that the
monodromy map
$\varphi$ can be chosen to be a symplectomorphism with a fixed point.

There is a related link surgery construction which starts with an oriented
$n$--component link
$L=\{K_1,\dots,K_n\}$ in $S^3$ and $n$ pairs $(X_i,T_i)$ of smoothly embedded
self-intersection $0$ tori in simply connected $4$--manifolds as above. Let
\[\alpha_L\co\pi_1(S^3\setminus L)\to \Z \]   denote the homomorphism
characterized by the property that it send the meridian $m_i$ of each component
$K_i$ to $1$. Let $N(L)$ be a tubular neighborhood of $L$.
Then if $\ell_i$ denotes the longitude of the component $K_i$, the curves
$\g_i=\ell_i + \alpha_L(\ell_i) m_i$ on
$\bd N(L)$ given by the $\alpha_L(\ell_i)$ framing of $K_i$ form the
boundary of a
Seifert surface for the link. In $S^1\x (S^3\setminus  N(L))$ let $T_{m_i}=S^1\x m_i$
and define the 4--manifold $X(X_1,\dots X_n;L)$ by
\[ X(X_1,\dots X_n;L)= (S^1\x (S^3\setminus  N(L))\cup\bigcup\limits_{i=1}^n (X_i\setminus 
(T_i\x
D^2))
\]
where $S^1\x\bd N(K_i)$ is
identified with $\bd N(T_i)$ so that for each $i$:
\[ [T_{m_i}]=[T_i], \ \ {\text{and}} \ \
[\g_i] = [{\text{pt}}\x\bd D^2]. \]

\begin{ther}{\rm\cite{KL4M}}\qua If each $T_i$ is homologically essential and
contained in a cusp
neighborhood in $X_i$ and if each $\pi_1(X\setminus T_i) =1$, then
$X(X_1,\dots X_n;L)$ is
simply connected and its Seiberg--Witten invariant is
\[\sw_{X(X_1,\dots
X_n;L)}=\DD_L(t_1,\dots,t_n)\cdot\prod_{j=1}^n\sw_{E(1)\#_{F=T_j}X_j}\]
where $t_j=\exp(2[T_j])$ and $\DD_L(t_1,\dots,t_n)$ is the symmetric
multivariable Alex\-ander
polynomial.
\end{ther}

\section{$2$--bridge knots}

Recall that 2--bridge knots, $K$, are classified by the double covers of
$S^3$ branched over
$K$, which are lens spaces.  Let $K(p/q)$ denote the 2--bridge knot whose
double branched cover
is the lens space $L(p,q)$. Here, $p$ is odd and $q$ is relatively prime to
$p$. Notice that
$L(p,q)\cong L(p,q-p)$; so we may assume at will that either $q$ is even or
odd.  We are first
interested in finding a pair of distinct fibered 2--bridge knots $K(p/q_i)$,
$i=1,2$ with the same Alexander polynomial.
Since 2--bridge knots are alternating, they are fibered if and only if their
Alexander
polynomials are monic \cite{BZ}. There is a simple combinatorial scheme for
calculating the
Alexander polynomial of a 2--bridge knot $K(p/q)$; it is described as
follows in \cite{S}.
Assume that $q$ is even and let $\mathbf{b}(p/q)=(b_1,\dots,b_n)$ where
$p/q$ is written as a
continued fraction:

\centerline{\unitlength 1cm
\begin{picture}(9,3)
\put (1,2.5){\LARGE $\frac{p}{q}$}
\put (1.35,2.5){ $=$}
\put (2,2.5){$2b_1 +1$}
\put (3,2.4){\line(1,0){.5}}
\put (2.8,2){$-2b_2 +1$}
\put (4.1,1.9){\line(1,0){.5}}
\put (4,1.5){$2b_3 +1$}
\put (5,1.4){\line(1,0){.5}}
\put (5.2,1.2){.}
\put (5.4,1){.}
\put (5.6,.8){.}
\put (5.7,.4){$+1$}
\put (6,.3){\line(1,0){.5}}
\put (6,-.1){$\pm 2b_n$}
\end{picture}}

There is then a Seifert surface for $K(p/q)$ whose corresponding Seifert
matrix is:
\[ V(p/q)=\begin{pmatrix} b_1 & 0 & 0 & 0 & 0 & \cdots \\ 1 & b_2 &1 & 0 &
0 &\cdots \\
        0 & 0 & b_3 & 0 & 0 &\cdots \\ 0 & 0 & 1 & b_4 & 1 &\cdots\\
        ..&..&..&..&..&\cdots
\end{pmatrix}\]
Thus the Alexander polynomial for $K(p/q)$ is
\[ \DD_{K(p/q)}(t) = \det(t\cdot V(p/q)-V(p/q)^{\text{tr}}).\]

\noindent Using this technique we calculate:

\begin{prop} The 2--bridge knots $K(105/64)$ and $K(105/76)$ share the
Alexander polynomial
\[ \DD(t)= t^4-5t^3+13t^2-21t+25-21t^{-1}+13t^{-2}-5t^{-3}+t^{-4}. \]
In particular, these knots are fibered.
\end{prop}
\proof The knots $K(105/64)$ and $K(105/76)$ correspond to the vectors
$$
{\bf{b}}(105/64)=(1,1,-1,-1,-1,-1,1,1)
$$
\vspace{-0.2in}
$$
{\bf{b}}(105/76)=(1,1,1,-1,-1,1,1,1).\eqno{\qed}
$$

\section{The examples}

Consider any pair of inequivalent fibered 2--bridge knots $K_i=K(p/q_i)$,
$i=1,2$, with the same
Alexander polynomial $\DD(t)$. Let $\tilde{K}_i=\pi_i^{-1}(K_i)$ denote the
branch knot in the
2--fold branched covering space $\pi_i\co L(p,q_i)\to S^3$, and let
$\tilde{m}_i=\pi_i^{-1}(m_i)$, with
$m_i$ the meridian of $K_i$.
Then
$M_{K_i}=S^1\x_{\varphi_i}\Sig$ with double  cover
   $\tilde{M}_{K_i} = S^1\x_{\varphi_i^2}\Sig$.

Let $X$ be the K3--surface and let $F$ denote a smooth torus of
self-intersection $0$ which is a
fiber of an elliptic fibration on $X$. Our examples are
\[X_{K_i} = X\#_{F=T_{\tilde{m}_i}}(S^1\x\tilde{M}_{K_i}). \]
The gluing is chosen so that the boundary of a normal disk to $F$ is
matched with the lift
$\tilde{\ell}_i$ of a longitude to $K_i$. A simple calculation and our
above discussion
implies that $X_{K_1}$ and $X_{K_2}$ are homeomorphic \cite{HK} and have
the same
Seiberg--Witten invariant:

\begin{thm}\label{SW} The manifolds $X_{K_i}$ are homeomorphic
symplectic rational
homology K3--surfaces
with fundamental groups $\pi_1(X_{K_i})=\Z_p$. Their
Seiberg--Witten invariants are
\[ \sw_{X_{K_i}}=\det(\varphi^2_{i,*}-\t^2 I)=\DD(\t)\cdot\DD(-\t) \]
where $\t=\exp([F])$.
\end{thm}

\section{Their universal covers}

The purpose of this final section is to prove our main theorem.

\begin{thm}\label{main}  $X_{K(105/64)}$ and
$X_{K(105/76)}$ are homeomorphic but not diffeomorphic
symplectic $4$--manifolds with the same
Seiberg--Witten invariant.
\end{thm}

Let $K_1=K(105/64)$ and $K_2=K(105/76)$. We have already shown that $X_{K_1}$
and $X_{K_2}$ are
homeomorphic symplectic $4$--manifolds with the same Seiberg--Witten invariant.
Suppose that $f\co X_{K_1}\to X_{K_2}$ is a diffeomorphism. It then satisfies
$f_*(\sw_{X_{K_1}}) = \sw_{X_{K_2}}$. Since these are both Laurent
polynomials in
the single variable $\t=\exp([F])$, and $[F]=[T_{\tilde{m}_i}]$ in $X_{K_i}$,
after appropriately orienting $T_{\tilde{m}_2}$, we must have
\[ f_*[T_{\tilde{m}_1}] = [T_{\tilde{m}_2}]. \]

We study the induced diffeomorphism $\f\co\hat{X}_{K_1}\to\hat{X}_{K_2}$ of
universal
covers. The universal cover $\hat{X}_{K_i}$ of $X_{K_i}$ is obtained as
follows.
Let $\vt_i\co S^3\to L(p,q_i)$ be the universal covering ($p=105$, $q_1=64$,
$q_2=76$) which induces the universal covering $\hat{\vt}_i\co\hat{X}_{K_i} \to
{X}_{K_i}$ , and let
$\hat{L}_i$ be the
$p$--component link
$\hat{L}_i=\vt_i^{-1}(\tilde{K}_i)$. The composition of the maps
$\varphi\circ\vt_i\co S^3\to S^3$ is a dihedral covering space branched over
$K_i$, and the link
$\hat{L}_i=\hat{L}(p/q_i)$ is classically known as the `dihedral covering
link' of $K(p/q_i)$. This is a symmetric link, and in fact, the deck
transformations
$\t_{i,k}$ of the cover $\vt_i\co S^3\to L(p,q_i)$ permute the link
components. The
collection of linking numbers of $\hat{L}_i$ (the dihedral linking numbers of
$K(p/q_i)$) classify the 2--bridge knots \cite{BZ}. The universal cover
$\hat{X}_{K_i}$
is obtained via the construction
$\hat{X}_{K_i}= X(X_1,\dots X_p;L_i)$ of \S2, where each
$(X_i,T_i)=(K3,F)$. Hence it
follows from \S2 that
$$
\sw_{\hat{X}_{K_i}} =
\DD_{\hat{L}_i}(t_{i,1},\dots,t_{i,p})\cdot\prod_{j=1}^p\sw_{E(1)\#_{F}K3}
=\hskip 1in\mbox{}$$\vspace{-0.2in}
$$\mbox{}\hskip 1.5in \DD_{\hat{L}_i}(t_{i,1},\dots,t_{i,p})\cdot
\prod_{j=1}^p(t_{i,j}^{1/2}-t_{i,j}^{-1/2})
$$
where $t_{i,j}=\exp([2T_{i,j}])$ and $T_{i,j}$ is the fiber $F$ in the $j$th
copy of
K3. Let $L_{i,1},\dots,L_{i,p}$ denote the components of the covering link
$\hat{L}_i$ in $S^3$, and let $m_{i,j}$ denote a meridian to $L_{i,j}$.
Then $[T_{i,j}]=[S^1\x m_{i,j}]$ in $H_2(\hat{X}_{K_i};\Z)$, and so
$\hat{\vt}_{i_*}[T_{i,j}] = [T_i]$.

Now we have $\f_*(\sw_{\hat{X}_{K_1}})=\sw_{\hat{X}_{K_2}}$ as elements of the integral
group ring of  $H_2(\hat{X}_{K_2};\Z)$. The formula given for
$\sw_{\hat{X}_{K_i}}$
shows that each basic class may be written in the form
$\b=\sum_{j=1}^p a_j[T_{i,j}]$.
Thus if $\b$ is a basic class of $\hat{X}_{K_1}$, then
\[
\f_*(\b) = \f_*(\sum\limits_{j=1}^p a_j[T_{1,j}]) =
\sum_{j=1}^p b_j[T_{2,j}]
\]
for some integers, $b_1,\dots,b_p$. But since $f_*[T_1]=[T_2]$ in
$H_2(X_{K_2};\Z)$ we have
$$
(\sum_{j=1}^p a_j)[T_2] = f_*(\sum_{j=1}^p a_j[T_1]) = f_*\hat{\vt}_{1*}(\b)
=\hat{\vt}_{2*}\f_*(\b)=\hskip 1in\mbox{}$$\vspace{-0.25in}
$$\mbox{}\hskip 1.5in\hat{\vt}_{2*}(\sum_{j=1}^p b_j[T_{2,j}])
= \sum_{j=1}^p b_j[T_2].
$$
Hence $\sum_{j=1}^p a_j = \sum_{j=1}^p b_j$.

Form the 1--variable Laurent polynomials
$P_i(t) = \DD_{\hat{L}_i}(t,\dots,t)\cdot (t^{1/2}-t^{-1/2})^p$
by equating all the variables $t_{i,j}$ in $\sw_{\hat{X}_{K_i}}$. The
coefficient
of a fixed term $t^k$ in $P_i(t)$ is
\[ \sum\{\ssw_{\hat{X}_{K_i}}(\sum_{j=1}^p a_j[T_{i,j}])\  \vert\
\sum_{j=1}^p a_j
= k\}. \]
Our argument above (and the invariance of the Seiberg--Witten invariant
under diffeomorphisms) shows that $\f_*$ takes $P_1(t)$ to $P_2(t)$; ie,
$P_1(t)=P_2(t)$ as Laurent polynomials.

The reduced Alexander polynomials $\DD_{\hat{L}_i}(t,\dots,t)$ have the form
\[\DD_{\hat{L_i}}(t,\dots,t)=(t^{1/2}-t^{-1/2})^{p-2}\cdot\nabla_{\hat
{L_i}}(t),\]where the polynomial $\nabla_{\hat{L_i}}(t
)$ is called the Hosokawa polynomial
\cite{H}.  Consider the matrix:
\[ \Lambda(p/q)=\begin{pmatrix} \s & \ve_1 & \cdots & \ve_{p-1}\\
       \ve_{p-1} & \s & \cdots & \ve_{p-2} \\
       : & : & &:\\  : & : & &:\\
       \ve_1 & \ve_2 & \cdots &\s
\end{pmatrix} \]
(Burde has shown that this is the linking matrix of $\hat{L}(p/q)$.)

It is a theorem of Hosokawa \cite{H} that
$\nabla_{\hat{L}(p/q)}(1)$ can be calculated as the determinant of any
$(p-1)$ by
$(p-1)$ minor $\Lambda^\prime(p/q)$ of $\Lambda(p/q)$. In particular, we
have the
following Mathematica calculations. (Note that
$K(105/64)=K(105/-41)$ and $K(105/76)=K(105/-29)$.)
\begin{gather*} \det(\Lambda^\prime(105/-41))/105 = 13^2\cdot 61^2\cdot
127^2\cdot
463^2\cdot
631^4\cdot 1358281^4 \\
\det(\Lambda^\prime(105/-29))/105 = 139^4\cdot 211^4\cdot 491^2\cdot
8761^2\cdot
10005451^4.
\end{gather*}

This means that $\nabla_{\hat{L_1}}(1) \ne \nabla_{\hat{L_2}}(1)$. However,
if we
let
$Q(t)=(t^{1/2}-t^{-1/2})^{2p-2}$, then
$P_i(t)=\nabla_{\hat{L_i}}(t)\cdot Q(t)$. For $|u-1|$ small enough,
$P_1(u)/Q(u) \ne P_2(u)/Q(u)$. Hence for $u\ne 1$ in this range, $P_1(u) \ne
P_2(u)$. This contradicts the existence of the diffeomorphism $f$ and completes
the proof of Theorem~\ref{main}.

\rk{Acknowledgements}The first author was partially supported NSF
Grant DMS9704927 and the second author by NSF Grant DMS9626330

\Addresses\recd

\end{document}